\documentclass[graybox]{svmult}
\usepackage[linesnumbered,ruled]{algorithm2e}
\usepackage{dsfont}
\usepackage{mathrsfs}  
\usepackage{enumitem}
\usepackage{bm}

\usepackage{type1cm}        
\usepackage{makeidx}         
\usepackage{graphicx}        
\usepackage{multicol}        
\usepackage[bottom]{footmisc}

\usepackage{newtxtext}       %
\usepackage{newtxmath}       

\makeindex             

\begin{document}

\title*{Construction of Grid Operators for Multilevel Solvers: a Neural Network Approach}
\author{Claudio Tomasi and Rolf Krause}
\institute{Claudio Tomasi and Rolf Krause \at Università della Svizzera Italiana,Via Buffi 13, CH-6904 Lugano, \email{claudio.tomasi@usi.ch, rolf.krause@usi.ch}}
%
%
\maketitle

\section{Introduction}
\label{Tomasi_mini_10_sec:1}

Multigrid (MG) methods are among the most successful strategies for solving linear systems arising from discretized elliptic equations. The main idea is to combine different levels of approximation in a multilevel hierarchy to compute the solution: it is possible to show that this algorithm is effective on the entire spectrum, thus leading to an optimal convergence property \cite{Tomasi_mini_10_braess2007finite, Tomasi_mini_10_briggs2000multigrid}.
Common to all these strategies is the need for the transfer of data or information between the different grids, or meshes.  Therefore, a crucial point for reaching fast convergence is the definition of transfer operators, but they are generally problem-dependent.
Except for the case of nested meshes, the computation of these operators is very expensive, and domain knowledge is always required.

The ever-increasing application of Machine Learning (ML) as support for methods in scientific computing makes it a natural solution to be employed in the definition of transfer operators, reducing the costs of their construction. In \cite{Tomasi_mini_10_greenfeld2019learning}, the learning of a mapping between PDEs and operators has been proposed. Another approach is presented in \cite{Tomasi_mini_10_katrutsa2020black}, where restriction and prolongation matrices are optimized while minimizing the spectral radius of the iteration matrix. As an alternative, the method proposed in \cite{Tomasi_mini_10_luz2020learning} uses Graph Neural Networks, for learning AMG prolongation operators, having classes of sparse matrices as input. 

In this paper, we propose a methodology based on Deep Neural Networks to define transfer operators based on the concept of $L^2$-projection. We take information from the domain to create several examples and to make our model learn from experience. Therefore, our focus is the construction of a suitable training set and a correct loss function definition to create a model that can be employed in MG solvers. The actual state of the method presents some limitations related to the mesh structure. An extension to a wider range of scenarios should be considered in future works.

\section{Problem Definition}
\label{Tomasi_mini_10_sec:2}

Let $\Omega \subset \mathbb{R}^n$ be a domain with Lipschitz boundary and let $H^1_0(\Omega)$ be the Sobolev space of one-time weakly differentiable functions on $\Omega$, with weak derivatives in $L^2(\Omega)$. We consider a multigrid method for the solution of the  following problem:

\begin{equation}
\mathrm{find} \ u \in V : a(u,v) = f(v) ~~ \forall v \in V
\label{Tomasi_mini_10_eq:1},
\end{equation}  
where $V\subset H_0^1(\Omega)$, $a : V \times V \rightarrow \mathbb{R}$ is a continuous symmetric elliptic bilinear form and $f: V \rightarrow \mathbb{R}$ is a continuous linear functional. Let $V_h \subset V$ be the associated finite elements space, where $\mathrm{dim}(V_h) = n_h$ and $h > 0$, and consider a conforming shape-regular triangulation $\mathcal{T}_h$. For a more rigorous explanation see e.g. \cite{Tomasi_mini_10_quarteroni2009numerical}.
Furthermore, let $I_H^h : \mathbb{R}^{n_H} \rightarrow  \mathbb{R}^{n_h}$ be a transfer operator which transfers information between $\mathcal{T}_h$ and a coarser grid $\mathcal{T}_H$, with $H > h$ and $n_H < n_h$. We denote with $A_h u_h = b_h$ the linear system arising from the finite element discretization  of (\ref{Tomasi_mini_10_eq:1}).

Let us consider a 2-grid correction scheme for solving the linear system. The extension to a general multigrid scenario is straightforward. To restrict or prolong information between coarse and fine grids, we apply $I_H^h$. Moreover, we define the coarse problem using the expression $A_H = (I_H^h)^T A_h~ I_H^h$. 
Hence, the definition of the transfer operator plays a central role in obtaining a fast convergence of the method.
In \cite{Tomasi_mini_10_krause2016parallel}, a general definition of transfer operators between meshes is discussed.
 Here, we focus on the $L^2$-projection as transfer operator. Let us call it $Q$:
\[	Q = M_h^{-1} B_h	,\]
where $M_h$ is the mass matrix related to the fine level (grid), and $B_h$ is a rectangular coupling operator matrix. The latter relates the two meshes, and it is computed through their intersection.
Since the inverse of $M_h$ is a dense matrix, the computation of $Q$ might become expensive. Therefore, we use the pseudo-$L^2$-projection, where we invert the lumped mass matrix instead of $M_h$.
For further reading refer to \cite{Tomasi_mini_10_dickopf2010multilevel, Tomasi_mini_10_dickopf2009pseudo, Tomasi_mini_10_dickopf2011study, Tomasi_mini_10_dickopf2014evaluating}.

\subsection{Neural Networks}
\label{Tomasi_mini_10_subsec:22}

ML algorithms are able to learn from data \cite{Tomasi_mini_10_goodfellow2016deep}; we refer to a single data object calling it example. An example is a collection of features together with a corresponding target. A feature is a property that has been measured from some object or event. The target is the correct response to the features, that the system should be able to reproduce. We represent an example as a couple $(\mathbf{x},\mathbf{y})$, where $\mathbf{x} \in \mathbb{R}^n$ is the feature set and $\mathbf{y} \in \mathbb{R}^m$ is the target.
ML can solve  different tasks, as classification, transcription, and so on. Our  focus is regression: we ask the model to predict numerical values given some inputs. In order to solve this task, the model is asked to output a function $f:\mathbb{R}^n \rightarrow \mathbb{R}^m$.
To evaluate the ML algorithm abilities, we define  a measure of its performance, called loss  function: for regression, we select the Mean Squared Error (MSE) indicator.

Neural Networks (NNs) belong to the class of supervised ML algorithms. They consist of layers of neurons, connected by weighted synapses. A NN defines a mapping $\mathbf{y} = f(\mathbf{x}; \bm\theta)$ and learns the values of the parameters $\bm \theta$, providing the best function approximation. More details can be found in \cite{Tomasi_mini_10_bishop2006pattern, Tomasi_mini_10_larochelle2009exploring}.

\subsection{Training Trasfer Operators}
\label{Tomasi_mini_10_subsec:23}

We aim to define a NN model to learn and then predict the transfer operator  $Q$.  Specifically, we do not learn directly $Q$, but the coupling operator $B_h$. Once the model is optimized, we employ it as a black box for solving linear systems of equations in an MG fashion.
We proceed by coarsening: we take $M_h$ on the fine level, and we extract the features in input to the NN. More details on the data extraction from $M_h$ are given in Section \ref{Tomasi_mini_10_subsec:31}. The model produces parts of $B_h$ that combined give rise to the full operator. We then retrieve the transfer operator $Q$, and we employ it in the MG algorithm. Furthermore, we use the predicted transfer operators to define the coarser mass matrix $M_H$ using the so-called \emph{Galerkin operator}, i.e.,
$M_H = Q^\top M_h Q. $
 \\We recursively apply this procedure to define coarser problems, giving rise to a multilevel hierarchy.
%

\section{Training Set}
\label{Tomasi_mini_10_sec:3}
In order to allow the NN to learn, we provide a large number of examples (or records). We need several distinct examples to be sure of avoiding overfitting, occurring when the model predictions correspond too closely or exactly to a particular set of data. Thus, we define classes of examples, and we choose a fixed  amount of records for each class. This allows us to create an unbiased training set without preferring some classes over others.
The definition of class is related to the mesh from which we extrapolate the records. We set a number of elements $N$: all the records coming from meshes with $N$ elements belong to class $C_N$.
\subsection{Records}
\label{Tomasi_mini_10_subsec:31}
Given a subset of pre-identified coarse nodes, we extract a record for each of them. Let $j$ be a coarse node. An example contains information related to \emph{patch($j$)}; here, \emph{patch($j$)} is the set consisting of node $j$ together with its neighbors. For each node $k \in$ \emph{patch($j$)}, we define features and target as the non-zero entries of the $k$th rows of  $M_h$ and  $B_h$, respectively.

 We consider different examples of 2-grid scenarios, where we associate to each fine mesh one coarse mesh, in order to approximate an actual function. Since we consider several examples for the same class $C_N$,  we need a strategy to avoid duplication inside the dataset. For this purpose, in each example, we consider the fine mesh and we move the nodes along the edges by a random quantity, proportional to the step size $h$. Therefore, we create different elements and consequently different records.
 We generate the examples in $C_N$, and we proceed to the next class by increasing $N$.
  Since NN models allow only fixed input and output dimensions, we define distinct models for 1D and 2D scenarios.

\subsection{One-Dimensional Model}
\label{Tomasi_mini_10_subsec:32}
The records related to one-dimensional meshes are extracted from scenarios obtained by coarsening: starting from a randomly generated fine mesh, we decide which nodes to keep for defining the coarse grid.  Here, \emph{patch($j$)}  consists only of $j$ itself, together with its left and right neighbors. For each coarse node, we take the information on \emph{patch($j$)}  from $M_h$ and $B_h$ following the strategy explained in Section \ref{Tomasi_mini_10_subsec:31}, to define each example.

\subsection{Two-Dimensional Model}
\label{Tomasi_mini_10_subsec:33}
Let us call \emph{patch-size} the number of nodes in a specific patch: given a node $j$, its patch-size is the cardinality of the set \emph{patch($j$)}.
In 2D, even in the same mesh, we can have nodes with different patch-sizes.  Hence, we start considering only a fixed triangulation, such that the nodes would have the same neighborhood pattern. We will focus on dealing with different patch-sizes later in the paper, referring to their treatment in the context of NNs.

A crucial point is to find a correct distribution of data in the training set, in terms of magnitude of the values. Since the NN should not prefer some examples - thus, some classes - over others, we need to define a correct and even filling of the training set. As a first approach, we relate the concept of class to the procedure of refinement. Mesh refinement is a strategy to increase the accuracy of the solution of a discretized problem. It works as an iterative procedure applied to the single elements of a mesh. Here we consider two different strategies: bisection, which halves each element, and mid-point refinement, which takes the mid-points of each edge and joins them to create new elements. 
When we refine, we deal with a new class of examples. Applying a training algorithm on these data results in a poor ability of approximation and a large prediction error, making a NN model unfit to work in a MG setting.  The refinement procedure makes the number of elements scale by a factor of 2 (bisection) or 4 (mid-point). In terms of domain of training examples, this means that the initial classes of records, i.e., $C_N$ with $N$ small, are close to each other, while their distance grows when $N$ increases. This turns out to produce an uneven training set, without a good balance in terms of data distribution. 
For this reason, we need a linear increase in the number of elements. If the classes of examples are evenly spaced in terms of domain, the network does not prefer some classes over others. 
Therefore, a second approach changes the definition of class, independent of the concept of refinement: we start from a number  $N$, and we create a mesh having exactly $N$ elements. Once we extract enough records, we proceed to the next class, increasing $N$ by a constant $K$, and create a new mesh with $N+K$ elements. For each class, we extract the same number of examples. Following this simple procedure, the resulting training set is effectively unbiased and with a good distribution of the examples. A learning algorithm applied to these data produces the expected good approximation. Therefore, a model trained on this dataset can be applied inside an MG scenario.

\section{Model Training}
\label{Tomasi_mini_10_sec:4}
A NN optimizes its parameters in order to reduce the prediction error. Employing MSE as loss function results in good predictions, but the model does not gain a good generalization property, i.e., the ability to perform well on previously unobserved inputs. 

Regularization helps us overcome this issue: it reduces the hypothesis space, allowing the NN to have a higher probability of choosing the most correct function. We introduce in the loss function some penalty terms related to the domain knowledge. These terms force constraints during the training phase in order to respect properties that the transfer operator must satisfy.

\subsection{Regularization}

During the construction of the training set, for each coarse node $j$, we extract patches of $M_h$ and parts of $B_h$. We use this information to ask the model to force some rules on the rows of the predicted coupling operator. 

We define the $j$th predicted and actual rows of the operator $Q$ as

\begin{equation}
Q_j^\text{pred} = \frac{1}{\sum M_j }B_j^\text{pred}, \;\;\;\;\;\;\;\;\; Q_j^\text{true} = \frac{1}{\sum M_j }B_j^\text{true},
\end{equation}
where $\Phi_j$ denotes the $j$th row of the operator $\Phi = \{M_h, B_h\}$.

We know that the predicted transfer operator should preserve constants (more details can be found in \cite[Section 3.2]{Tomasi_mini_10_dickopf2010multilevel}). Hence, we consider the following penalty terms to specialize our loss function:
\begin{align}
\| Q_j^\text{pred} \cdot \mathds{1}_H -  \mathds{1}_h  \|_{_2}, &&  \| Q_j^\text{pred} - Q_j^\text{true}  \|_{_2},
\end{align}
where $\|\cdot\|_{_2}$ denotes the Euclidean norm,  $\mathds{1}_H$ and  $\mathds{1}_h$ the all-ones vectors of dimensions $n_H$ and  $n_h$, respectively. We then define, for all the nodes $k \in \;$\emph{patch($j$)}
\begin{equation}
p_k = \frac{1}{\alpha}\| Q_k^\text{pred} \cdot \mathds{1}_H -  \mathds{1}_h  \|_{_2} \,+\,   \frac{1}{\beta}\| Q_k^\text{pred} - Q_k^\text{true}  \|_{_2},
\end{equation}
where $0 < \alpha, \beta < 1$.

Therefore, we define the loss function as
\begin{equation}
 \mathscr{L}( y_\text{true}, y_\text{pred}) = MSE( y_\text{true}, y_\text{pred})  + \sum_k p_k .
\end{equation}
Adopting the latter during the training phase, in addition to minimize the simple distance between target and prediction, we aim to respect the above properties related to the transfer operator.
\subsection{Model Details}
We use a classic splitting for our dataset: $20\%$ for test and $80\%$ for training, where the latter is divided again in $20\%$ for the validation set and the remaining for the training phase. For a preliminary examination of the method, we used around $500.000$ examples. For both one- and two-dimensional models, we adopt Adam as optimizer. Regarding the architecture, we report here only the structure for the 2D scenario used as initial test: we need at least $20$ hidden layers, where for each of them we use $800$ neurons, for a total of  $12$ million parameters. Through further investigations and several tests, the NN complexity can be improved, giving rise to less expensive computations. We initialize our weights using a normal distribution, using the methods provided by Tensorflow. In the context of the NN definition, more extensive works should be devoted to study the sensitivity of the predictions while changing the NN parameters.

\section{Numerical Results}
\label{Tomasi_mini_10_sec:5}
We test our NNs for both prediction accuracy and their application in an MG setting. We compare our method  with the Semi-Geometric multigrid (SGMG) method (see \cite[Chapter 3]{Tomasi_mini_10_dickopf2010multilevel}), which adopts the $L^2$-projection computed through intersections between meshes. In addition to the convergence, we consider the difference in the time spent to assemble the transfer operator. For our method, we take into account patches extraction, predictions, assembly of $B$ and the computation of the operator. Regarding the computation of the actual $L^2$-projection we consider the time spent for intersecting fine and coarse mesh, triangulation for each intersecting polygon and numerical integration.
We test the method on one-dimensional examples, and the results are good as expected, considering two or more levels. Our method converges with the same number of iterations as the SGMG method. Comparing the CPU time spent in creating the transfer operators, we see that the predicted one is assembled faster  than the other since it depends only on the problem dimensions.

During the test on two-dimensional settings we need to deal with different patch-sizes, as described in Section \ref{Tomasi_mini_10_subsec:33}. Even if we consider a simple regular mesh, the nodes near the boundaries have fewer neighbors than the internal nodes. A preliminary solution requires the mesh to be extended, to make all the nodes have the same patch-size. Virtually, we add neighbors to those nodes having a smaller patch in the given mesh. Using this expedient, the method works, and we can test the convergence against the SGMG method. 
Extending the mesh shows to be useful for an initial application of MG, but it is very expensive in terms of computations. Increasing the degrees of freedom (dofs), we would have more and more virtual nodes to add and heavier computations to carry out. Therefore, we consider different NNs, each of them defined and optimized for a specific patch-size. 

\begin{figure}[t]
	\centering
	\includegraphics[scale=.26]{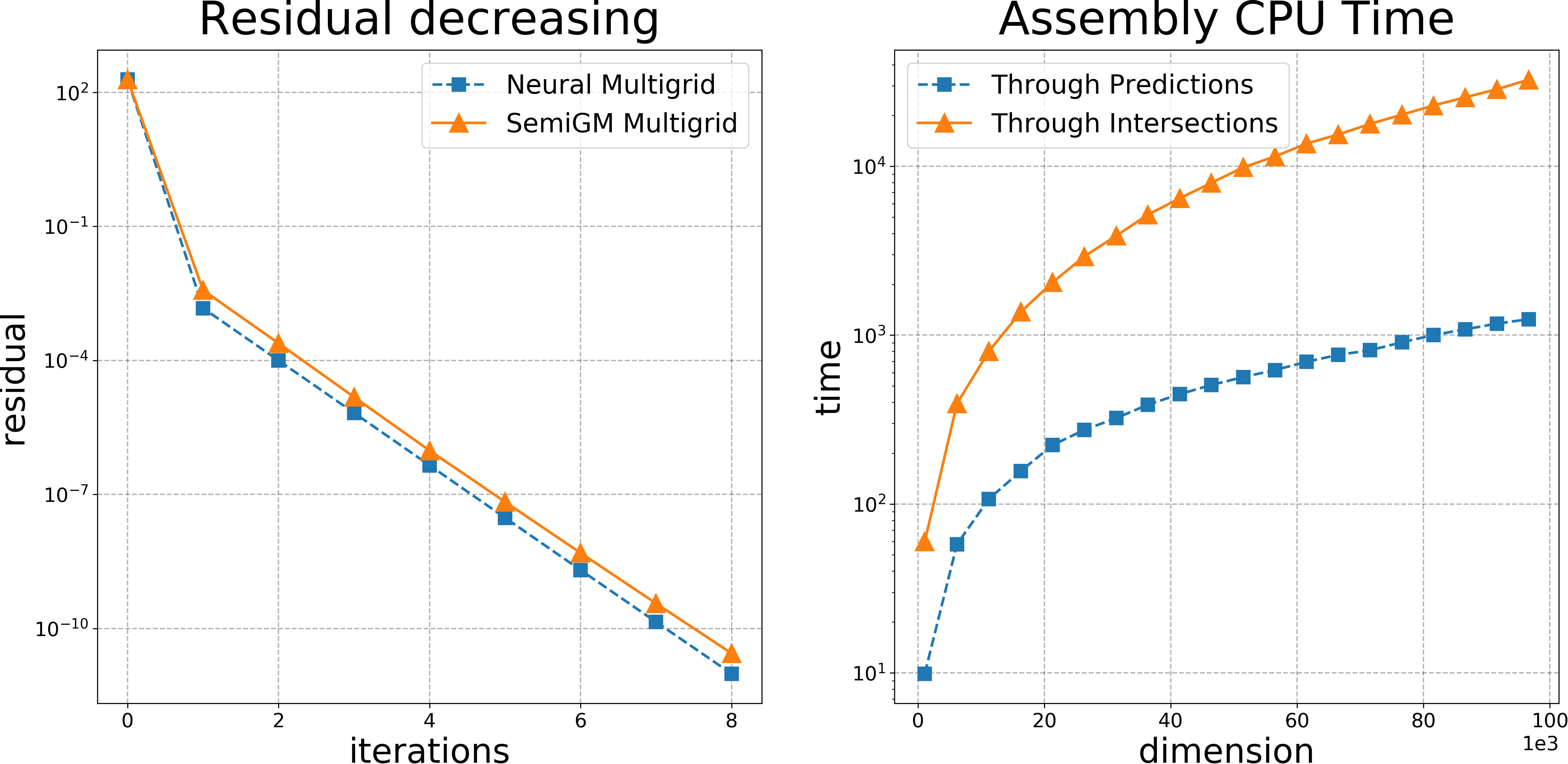}
	\caption{Convergence of 2D Neural MG against SGMG on an example of 100.000 dofs (left). Comparison of CPU time between the two methods, increasing the dofs (right).}
	\label{Tomasi_mini_10_fig:2D}       
\end{figure}

Fig. \ref{Tomasi_mini_10_fig:2D} shows the performance of the method on a two-dimensional scenario: in the left picture, we compare the convergence of our method  against the convergence of Semi-Geometric MG; in the right picture we compare the CPU time spent in both methods.

\section{Conclusion}
\label{Tomasi_mini_10_sec:6}
This work presented the study and definition of a methodology to construct NNs to predict transfer operators for MG solvers. Starting from a one-dimensional case, we built an unbiased training set allowing the optimization of a model, which brought very good results in an MG context. Reproducing the same methodology, we approached the two-dimensional setting, which gave us the chance to better define a training set for this kind of methods. Furthermore, we could test our method using different input-sized neural networks, resulting in fast convergence and bringing a great speedup in the computation of the transfer operator.  
The same procedure can be employed for constructing models to deal with a general $N$-dimensional scenario. 
Given the limitations of this method at its current state, further investigations should be devoted to overcome the necessity of having multiple NNs modeled on different patch-sizes, in order to define a general strategy for solving arbitrary problems.
Future works should extend the method to deal with a wider class of triangulation, and for applications in other Multilevel scenarios.

\bibliography{Tomasi_mini_10.bib}
\bibliographystyle{abbrv}
\end{document}